\documentclass{article}
\usepackage{amsfonts}
\usepackage{amsmath}

\newtheorem{theorem}{Theorem}

\input{tcilatex}
\begin{document}

\title{Killing vector fields of the metric $II+III$ on Tangent Bundle}
\author{Melek Aras\thanks{%
Department of Mathematics,Faculty of Arts and Sciences, Giresun University,
28049, Turkey e-mail:melekaras25@hotmail.com;melek.aras@giresun.edu.tr }}
\maketitle

\begin{abstract}
The main purpose of the paper is to investigate Killing vector fields on the
tangent bundle $T\left( M_{n}\right) $ of the Riemannian manifold with
respect to the Levi-Civita connection of the metric $II+III$ .


2000 AMS Classification 53C07, 53C25
\end{abstract}

\begin{description}
\item 
\begin{description}
\item[1.Introduction] 
\end{description}
\end{description}

Let $M_{n}$ be an $n-$dimensional differentiable Riemannian manifold of
class $C^{\infty }$ with $g$ . Then the set $T\left( M_{n}\right) $ is
tangent bundle over the manifold $M_{n}\left( see\cite{2}\right) $. We
denote by $\Im _{q}^{p}\left( M_{n}\right) $ the set of all\ tensor fields
of type $\left( p,q\right) $ in $M_{n}$ and by $\pi :T\left( M_{n}\right)
\rightarrow M_{n}$ the naturel \ bundle structure of $T\left( M_{n}\right) $
over $M_{n}$. For $U\subset M_{n},\left( x^{i},x^{i^{\prime }}\right)
,i=1,...,n$ and $i^{\prime }=n+1,...,2n$ are local coordinates in a
neighborhood $\pi ^{-1}\left( U\right) \subset T\left( M_{n}\right) \left(
see\cite{6}\right) $.

\ \ Let $M_{n}$ be a Riemannian manifold with metric $g$ whose components in
a coordinate neighborhood $U$ are $g_{ji}$. In the neighborhood $\pi
^{-1}\left( U\right) $ of $T\left( M_{n}\right) $, $U$ being a neighborhood
of $M_{n}$, we put

\begin{equation*}
\delta y^{h}=dy^{h}+\Gamma _{i}^{h}dx^{i}
\end{equation*}%
with respect to the induced coordinates $\left( x^{h},y^{h}\right) $ in $\pi
^{-1}\left( U\right) \subset T\left( M_{n}\right) $, where $\Gamma
_{i}^{h}=y^{j}$ $\Gamma _{ji}^{h}$ \cite{6}.

Let $g$ be a Riemannian metric of $M_{n}$ with components $g_{ji}$, then we
see that

\begin{equation*}
\ \ II+III:2g_{ji}dx^{j}\delta y^{i}+g_{ji}\delta y^{j}\delta y^{i}
\end{equation*}%
\ is non-singular and consequently can be regarded as Riemannian or
pseudo-Riemannian metric on the tangent bundle $T\left( M_{n}\right) $ over $%
M_{n}$ $\left( see\cite{6}\right) $.

The metric $II+III$ has components\cite{3}

\begin{equation}
\ \ II+III:\left( \widetilde{g}_{CB}\right) =\left( 
\begin{array}{cc}
0 & g_{ji} \\ 
g_{ji} & g_{ji}%
\end{array}%
\right)  \tag{1}
\end{equation}%
and consequently its contravariant components

\begin{equation}
\left( \widetilde{g}^{CB}\right) =\left( 
\begin{array}{cc}
-g^{ji} & g^{ji} \\ 
g^{ji} & 0%
\end{array}%
\right)  \tag{2}
\end{equation}%
with respect to the adapted frame on $T\left( M_{n}\right) $.

The frame components of Levi-Civita connection of lift metric $\widetilde{g}$
are as follows\cite{3}:

\begin{equation}
\left\{ 
\begin{array}{c}
\widetilde{\Gamma }_{ji}^{h}=\Gamma _{ji}^{h}-\frac{1}{2}y^{b}\left(
R_{bji}^{h}+R_{bij}^{h}\right) ,\ \ \ \widetilde{\Gamma }_{ji}^{\overline{h}%
}=y^{b}R_{bji}^{h},\ \ \widetilde{\Gamma }_{\overline{ji}}^{\overline{h}%
}=0,\ \ \widetilde{\Gamma }_{ji}^{h}=0 \\ 
\widetilde{\Gamma }_{j\overline{i}}^{\overline{h}}=\Gamma _{ji}^{h}+\frac{1}{%
2}y^{b}R_{bij}^{h},\widetilde{\Gamma }_{j\overline{i}}^{h}=-\frac{1}{2}%
y^{b}R_{bij}^{h},\widetilde{\Gamma }_{\overline{j}i}^{\overline{h}}=\frac{1}{%
2}y^{b}R_{bji}^{h},\widetilde{\Gamma }_{\overline{j}i}^{h}=-\frac{1}{2}%
y^{b}R_{bji}^{h}%
\end{array}%
\right.  \tag{3}
\end{equation}%
where $\ \Gamma _{ji}^{h}$ denote the Christoffel symbols constructed with $%
g_{ji}$ on $M_{n}$.

Let $\widetilde{X}$ be a vector field in $T\left( M_{n}\right) $ and $\left( 
\widetilde{X}\right) =\left( 
\begin{array}{c}
\widetilde{X}^{h} \\ 
\widetilde{X}^{\overline{h}}%
\end{array}%
\right) $ its components with respect to the adapted frame. Then covariant
derivative $\widetilde{\nabla }\widetilde{X}$ has components

\begin{equation}
\widetilde{\nabla }_{\gamma }\widetilde{X}^{\alpha }=D_{\gamma }\widetilde{X}%
^{\alpha }+\widetilde{\Gamma }_{\gamma \beta }^{\alpha }\widetilde{X}^{\beta
},  \tag{4}
\end{equation}%
$\widetilde{\Gamma }_{\gamma \beta }^{\alpha }$ being given by $\left(
3\right) $, with respect to the adapted frame \cite{6}.

Consider a vector field $X$ in $M_{n}.$Then its vertical lift $^{V}X$,
complete lift $^{C}X$ and horizontal lift $^{H}X$ have respectively
components\cite{6}

\begin{equation}
\left( ^{^{\prime }}X^{A}\right) =\left( 
\begin{array}{c}
0 \\ 
X^{h}%
\end{array}%
\right) ,\left( \widetilde{X}^{A}\right) =\left( 
\begin{array}{c}
X^{h} \\ 
\partial X^{h}%
\end{array}%
\right) ,\left( \overline{X}^{A}\right) =\left( 
\begin{array}{c}
X^{h} \\ 
-\Gamma _{i}^{h}X^{i}%
\end{array}%
\right)  \tag{5}
\end{equation}%
with respect to with respect to the induced coordinates in $T\left(
M_{n}\right) .$ Then their components\ \cite{6}

\begin{equation*}
^{\prime }X^{\alpha }=A_{A}^{\alpha \prime }X^{A},\ \widetilde{X}^{\alpha
}=A_{A}^{\alpha }\widetilde{X}^{A}\ ,\ \ \ \overline{X}^{\alpha
}=A_{A}^{\alpha }\overline{X}^{A}
\end{equation*}%
with respect to the adapted frame are given respectively by

\begin{equation}
\left( ^{^{\prime }}X^{\alpha }\right) =\left( 
\begin{array}{c}
0 \\ 
X^{h}%
\end{array}%
\right) ,\left( \widetilde{X}^{\alpha }\right) =\left( 
\begin{array}{c}
X^{h} \\ 
\nabla X^{h}%
\end{array}%
\right) ,\left( \overline{X}^{\alpha }\right) =\left( 
\begin{array}{c}
X^{h} \\ 
0%
\end{array}%
\right)  \tag{6}
\end{equation}%
where $\nabla X^{h}=y^{i}\nabla _{i}X^{h}$(see Sasaki \cite{1}).

\begin{description}
\item[2. Killing vector fields of the metric II+III] 
\end{description}

Now let us consider the covariant derivatives of vertical, complete and
horizontal lifts of the vector field $X$ in $M_{n}.$

The covariant derivative $\widetilde{\nabla }^{V}X$ has components

\begin{equation}
\left( \widetilde{\nabla }_{\beta }^{V}X^{\alpha }\right) =\left( 
\begin{array}{cc}
-\frac{1}{2}y^{s}R_{sij}^{h}X^{j} & 0 \\ 
\nabla _{i}X^{h}+\frac{1}{2}y^{s}R_{sij}^{h}X^{j} & 0%
\end{array}%
\right)  \tag{7}
\end{equation}%
with respect to the adapted frame,because of $\left( 3\right) $ and $\left(
6\right) $. Thus we have

\begin{theorem}
The vertical lift of a vector field in $M_{n}$ to $T\left( M_{n}\right) $
with the metric \ $II+III$ is parallel if and only if the given vector field
in $M_{n}$ \ is parallel(see\cite{5}).
\end{theorem}

\bigskip The covariant derivative $\widetilde{\nabla }^{C}X$ has components

\begin{center}
$\widetilde{\nabla }_{\beta }^{C}X^{\alpha }=\left( 
\begin{array}{cc}
\widetilde{\nabla }_{i}^{C}X^{h} & \widetilde{\nabla }_{\overline{i}%
}^{C}X^{h} \\ 
\widetilde{\nabla }_{i}^{C}X^{\overline{h}} & \widetilde{\nabla }_{\overline{%
i}}^{C}X^{\overline{h}}%
\end{array}%
\right) $
\end{center}

\begin{equation}
\left\{ 
\begin{array}{l}
\widetilde{\nabla }_{i}^{C}X^{h}=%
\begin{array}{cc}
\nabla _{i}X^{h} & -\frac{1}{2}y^{s}\left( R_{sji}^{h}+R_{sij}^{h}\right)
X^{j}-\frac{1}{2}y^{s}R_{sij}^{h}\nabla _{l}X^{j}y^{l}%
\end{array}
\\ 
\widetilde{\nabla }_{\overline{i}}^{C}X^{h}=-\frac{1}{2}y^{s}R_{sji}^{h}X^{j}
\\ 
\widetilde{\nabla }_{i}^{C}X^{\overline{h}}=\nabla _{i}\nabla
_{l}X^{h}y^{l}+y^{s}R_{sij}^{h}X^{j}+\frac{1}{2}y^{s}R_{sij}^{h}\nabla X^{j}
\\ 
\widetilde{\nabla }_{\overline{i}}^{C}X^{\overline{h}}=\nabla _{i}X^{h}+%
\frac{1}{2}y^{s}R_{sji}^{h}X^{j} \\ 
\end{array}%
\right.  \tag{8}
\end{equation}%
with respect to the adapted frame,because of $\left( 3\right) $ and $\left(
6\right) $. Thus we have

\begin{theorem}
The complete lift of a vector field in $M_{n}$ to $T\left( M_{n}\right) $
with the metric \ $II+III$ is parallel if and only if the given vector field
in $M_{n}$ \ is parallel(see\cite{5}).
\end{theorem}

The covariant derivative $\widetilde{\nabla }^{H}X$ has components

\begin{equation}
\left( \widetilde{\nabla }_{\beta }^{H}X^{\alpha }\right) =\left( 
\begin{array}{cc}
\nabla _{i}X^{h}-\frac{1}{2}y^{s}\left( R_{sji}^{h}+R_{sij}^{h}\right) X^{j}
& -\frac{1}{2}y^{s}R_{sji}^{h}X^{j} \\ 
y^{s}R_{sji}^{h}X^{j} & \frac{1}{2}y^{s}R_{sji}^{h}X^{j}%
\end{array}%
\right)  \tag{9}
\end{equation}%
with respect to the adapted frame,because of $\left( 3\right) $ and $\left(
6\right) $.Thus we have

\begin{theorem}
The horizontal lift of a vector field in $M_{n}$ to $T\left( M_{n}\right) $
with the metric \ $II+III$ are parallel if and only if the given vector
field in $M_{n}$ \ is parallel(see\cite{5}).
\end{theorem}

Given a vector field $\widetilde{X}$ in $T\left( M_{n}\right) $, the \ $%
1-form$ $\widetilde{\omega }$ defined by $\widetilde{\omega }\left( 
\widetilde{Y}\right) =\widetilde{g}\left( \widetilde{X},\widetilde{Y}\right)
,$ $\widetilde{Y}$ being an arbitrary element of $\Im _{0}^{1}\left( T\left(
M_{n}\right) \right) ,$ is called the covector field associated with $%
\widetilde{X}$ \ and denoted by $\widetilde{X}$ $^{\ast }$. If $\ \widetilde{%
X}$ \ has local components $\widetilde{X}$ $^{A}$, then the associated
covector field with $\widetilde{X}$ $^{\ast }$ of $\ \widetilde{X}$ \ has
local components $\widetilde{X}_{C}=\widetilde{g}_{CB}\widetilde{X}^{B}\cite%
{6}.$

The vertical , complete and horizontal lifts of a vector field $X$ in $M_{n}$
with components $X^{h}$ are given by $\left( 6\right) $.Then using $\left(
2\right) ,$we see that covector fields respectively associated with $%
^{V}X,^{C}X$ and $^{H}X$ have in $T\left( M_{n}\right) $ with the metric $%
II+III$ \ components

\begin{equation}
\left( ^{V}X_{\beta }\right) =\left( X_{i},X_{i}\right) \text{, \ \ }\left(
^{C}X_{\beta }\right) =\left( \nabla X_{i},X_{i}+\nabla X_{i}\right) \text{,
\ \ }\left( ^{H}X_{\beta }\right) =\left( 0,X_{i}\right)  \tag{10}
\end{equation}%
with respect to the adapted frame, where $X_{i}=g_{ih}X^{h}$ are components
of the covector field $X^{\ast }$ associated with $X$ . Thus the rotations
of $^{V}X$, $^{C}X$ and $^{H}X$ have respectively components

\begin{equation}
\left( \widetilde{\nabla}_{\beta}^{V}X_{\alpha}-\widetilde{\nabla}_{\alpha
}^{V}X_{\beta}\right) =\left( 
\begin{array}{cc}
\nabla_{i}X_{j}-\nabla_{j}X_{i}+\left( R_{sji}^{h}-R_{sij}^{h}\right)
y^{s}X_{h} & \nabla_{i}X_{j}-\nabla_{j}X_{i} \\ 
0 & 0%
\end{array}
\right)  \tag{11}
\end{equation}

\begin{center}
$\left( \widetilde{\nabla}_{\beta}^{C}X_{\alpha}-\widetilde{\nabla}_{\alpha
}^{C}X_{\beta}\right) =\left( 
\begin{array}{cc}
\widetilde{\nabla}_{i}^{C}X_{j}-\widetilde{\nabla}_{j}^{C}X_{i} & \widetilde{%
\nabla}_{i}^{C}X_{\overline{j}}-\widetilde{\nabla}_{\overline{j}}^{C}X_{i}
\\ 
\widetilde{\nabla}_{\overline{i}}^{C}X_{j}-\widetilde{\nabla}_{j}^{C}X_{%
\overline{i}} & \widetilde{\nabla}_{\overline{i}}^{C}X_{\overline{j}}-%
\widetilde{\nabla}_{\overline{j}}^{C}X_{\overline{i}}%
\end{array}
\right) $
\end{center}

\begin{equation}
\left\{ 
\begin{array}{l}
\widetilde{\nabla }_{i}^{C}X_{j}-\widetilde{\nabla }_{j}^{C}X_{i}=\left(
\nabla _{i}\nabla _{l}-\nabla _{l}\nabla _{i}\right) y^{l}X_{j}+\left(
R_{sji}^{h}-R_{sij}^{h}\right) y^{s}\nabla X_{h}+\left(
R_{sij}^{h}-R_{sji}^{h}\right) y^{s}X_{h} \\ 
\widetilde{\nabla }_{i}^{C}X_{\overline{j}}-\widetilde{\nabla }_{\overline{j}%
}^{C}X_{i}=\left( \nabla _{i}X_{j}-\nabla _{j}X_{i}\right) +\left( \nabla
_{i}\nabla _{l}-\nabla _{l}\nabla _{i}\right) y^{l}X_{j} \\ 
\widetilde{\nabla }_{\overline{i}}^{C}X_{j}-\widetilde{\nabla }_{j}^{C}X_{%
\overline{i}}=-\frac{1}{2}y^{s}\left( R_{sij}^{h}-R_{sji}^{h}\right) X_{h}
\\ 
\widetilde{\nabla }_{\overline{i}}^{C}X_{\overline{j}}-\widetilde{\nabla }_{%
\overline{j}}^{C}X_{\overline{i}}=0 \\ 
\end{array}%
\right.  \tag{12}
\end{equation}

\begin{center}
$\left( \widetilde{\nabla}_{\beta}^{H}X_{\alpha}-\widetilde{\nabla}_{\alpha
}^{H}X_{\beta}\right) =\left( 
\begin{array}{cc}
\widetilde{\nabla}_{i}^{H}X_{j}-\widetilde{\nabla}_{j}^{H}X_{i} & \widetilde{%
\nabla}_{i}^{H}X_{\overline{j}}-\widetilde{\nabla}_{\overline{j}}^{H}X_{i}
\\ 
\widetilde{\nabla}_{\overline{i}}^{H}X_{j}-\widetilde{\nabla}_{j}^{H}X_{%
\overline{i}} & \widetilde{\nabla}_{\overline{i}}^{H}X_{\overline{j}}-%
\widetilde{\nabla}_{\overline{j}}^{H}X_{\overline{i}}%
\end{array}
\right) $
\end{center}

\begin{equation}
\left( \widetilde{\nabla}_{\beta}^{H}X_{\alpha}-\widetilde{\nabla}_{\alpha
}^{H}X_{\beta}\right) =\left\{ 
\begin{array}{l}
\widetilde{\nabla}_{i}^{H}X_{j}-\widetilde{\nabla}_{j}^{H}X_{i}=\left(
R_{sij}^{h}-R_{sji}^{h}\right) y^{s}X_{h} \\ 
\widetilde{\nabla}_{i}^{H}X_{\overline{j}}-_{\overline{j}}^{H}X_{i}=\left(
\nabla_{i}X_{j}-\nabla_{j}X_{i}\right) +\frac{1}{2}\left(
R_{sji}^{h}-R_{sij}^{h}\right) y^{s}X_{h} \\ 
\widetilde{\nabla}_{\overline{i}}^{H}X_{j}-\widetilde{\nabla}_{j}^{H}X_{%
\overline{i}}=\frac{1}{2}\left( R_{sij}^{h}-R_{sji}^{h}\right) y^{s}X_{h} \\ 
\widetilde{\nabla}_{\overline{i}}^{H}X_{\overline{j}}-\widetilde{\nabla }_{%
\overline{j}}^{H}X_{\overline{i}}=0 \\ 
\end{array}
\right.  \tag{13}
\end{equation}

with respect to the adapted frame.

From $\left( 12\right) $, we see that

\begin{equation}
\nabla_{i}X_{j}-\nabla_{j}X_{i}=0,\text{ \ \ \ }\nabla_{i}\nabla_{l}X_{j}=0%
\text{,}  \tag{14}
\end{equation}

if the complete lift of $X^{\ast }$ is closed in $T\left( M_{n}\right) $.
Further, if the conditions $\left( 14\right) $ are satisfied, we easily
deduce that $\left( R_{sji}^{h}-R_{sij}^{h}\right) y^{s}\nabla X_{h}=0$ and $%
\left( R_{sij}^{h}-R_{sji}^{h}\right) y^{s}X_{h}=0,$ so that the complete
lift of $X^{\ast }$ is closed in $T\left( M_{n}\right) $ if $X^{\ast }$ is
closed and the second covariant derivative of $X$ \ vanishes in $M_{n}$.

Thus we have

\begin{theorem}
Necessary and sufficient condition in order that $\left( a\right) $ vertical,%
$\left( b\right) $ complete, $\left( c\right) $ horizontal lifts to $T\left(
M_{n}\right) $ with the metric $II+III$ , of a vector field $X$ in $M_{n}$
be harmonic in $T\left( M_{n}\right) $ are respectively that the given
vector field $X$ in $M_{n}$ $\left( a\right) $\ parallel; $\left( b\right) $
harmonic and having vanishing second covariant derivative ; and $\left(
c\right) $ harmonic\cite{5}.
\end{theorem}

A vector field $X\epsilon \Im _{0}^{1}$ $\left( M_{n}\right) $ is said to be
a Killing vector field of a Riemannian manifold with metric $g$, if $%
\tciLaplace _{X}g=0\cite{4}$. In terms of components $g_{ji}$ of $g$, $X$ is
a Killing vector field if only if

\begin{equation*}
\tciLaplace _{X}g=X^{\alpha }\nabla _{\alpha }g_{ji}+g_{\alpha i}\nabla
_{j}X^{\alpha }+g_{j\alpha }\nabla _{i}X^{\alpha }=\nabla _{j}X_{i}+\nabla
_{i}X_{j}=0,
\end{equation*}%
$X^{\alpha }$ being components of $X$, where $\nabla $ is the Riemannian
connection of the metric $g\cite{4}$.

The Lie derivatives of the metric $II+III$ with respect to $^{V}X$, $^{C}X$
and $^{H}X$ have respectively components

\QTP{Body Math}
\begin{equation}
\left( \widetilde{\nabla }_{\beta }^{V}X_{\alpha }+\widetilde{\nabla }%
_{\alpha }^{V}X_{\beta }\right) =\left( 
\begin{array}{cc}
\nabla _{i}X_{j}+\nabla _{j}X_{i} & \nabla _{i}X_{j}+\nabla _{j}X_{i} \\ 
0 & 0%
\end{array}%
\right)  \tag{11}
\end{equation}

\begin{center}
$\left( \widetilde{\nabla }_{\beta }^{C}X_{\alpha }+\widetilde{\nabla }%
_{\alpha }^{C}X_{\beta }\right) =\left( 
\begin{array}{cc}
\widetilde{\nabla }_{i}^{C}X_{j}+\widetilde{\nabla }_{j}^{C}X_{i} & 
\widetilde{\nabla }_{i}^{C}X_{\overline{j}}+\widetilde{\nabla }_{\overline{j}%
}^{C}X_{i} \\ 
\widetilde{\nabla }_{\overline{i}}^{C}X_{j}+\widetilde{\nabla }_{j}^{C}X_{%
\overline{i}} & \widetilde{\nabla }_{\overline{i}}^{C}X_{\overline{j}}+%
\widetilde{\nabla }_{\overline{j}}^{C}X_{\overline{i}}%
\end{array}%
\right) $

\begin{equation}
\left\{ 
\begin{array}{l}
\widetilde{\nabla }_{i}^{C}X_{j}+\widetilde{\nabla }_{j}^{C}X_{i}=\nabla
_{i}\nabla _{l}X_{j}y^{l}+\left( R_{hsij}+R_{hsji}\right) X^{h}y^{s} \\ 
\text{ \ }\widetilde{\nabla }_{i}^{C}X_{\overline{j}}+\widetilde{\nabla }_{%
\overline{j}}^{C}X_{i}=\left( \nabla _{i}X_{j}+\nabla _{j}X_{i}\right)
+\left( \nabla _{i}\nabla _{l}X_{j}+\nabla _{j}\nabla _{l}X_{i}\right) y^{l}
\\ 
\widetilde{\nabla }_{\overline{i}}^{C}X_{j}+\widetilde{\nabla }_{j}^{C}X_{%
\overline{i}}=-\frac{1}{2}y^{s}\left( R_{sij}^{h}+R_{sji}^{h}\right) X_{h}
\\ 
\widetilde{\nabla }_{\overline{i}}^{C}X_{\overline{j}}+\widetilde{\nabla }_{%
\overline{j}}^{C}X_{\overline{i}}=0 \\ 
\end{array}%
\right.  \tag{12}
\end{equation}

$\left( \widetilde{\nabla }_{\beta }^{H}X_{\alpha }+\widetilde{\nabla }%
_{\alpha }^{H}X_{\beta }\right) =\left( 
\begin{array}{cc}
\widetilde{\nabla }_{i}^{H}X_{j}+\widetilde{\nabla }_{j}^{H}X_{i} & 
\widetilde{\nabla }_{i}^{H}X_{\overline{j}}+\widetilde{\nabla }_{\overline{j}%
}^{H}X_{i} \\ 
\widetilde{\nabla }_{\overline{i}}^{H}X_{j}+\widetilde{\nabla }_{j}^{H}X_{%
\overline{i}} & \widetilde{\nabla }_{\overline{i}}^{H}X_{\overline{j}}+%
\widetilde{\nabla }_{\overline{j}}^{H}X_{\overline{i}}%
\end{array}%
\right) $
\end{center}

\begin{equation}
\left\{ 
\begin{array}{l}
\widetilde{\nabla }_{i}^{H}X_{j}+\widetilde{\nabla }_{j}^{H}X_{i}=y^{s}%
\left( R_{sij}^{h}+R_{sji}^{h}\right) X_{h} \\ 
\text{ \ }\widetilde{\nabla }_{i}^{H}X_{\overline{j}}+\widetilde{\nabla }_{%
\overline{j}}^{H}X_{i}=\nabla _{i}X_{j}+\nabla _{j}X_{i}+\frac{1}{2}\left(
R_{sji}^{h}+R_{sij}^{h}\right) y^{s}X_{h} \\ 
\widetilde{\nabla }_{\overline{i}}^{H}X_{j}+\widetilde{\nabla }_{j}^{H}X_{%
\overline{i}}=-\frac{1}{2}y^{s}\left( R_{sij}^{h}+R_{sji}^{h}\right) X_{h}
\\ 
\widetilde{\nabla }_{\overline{i}}^{H}X_{\overline{j}}+\widetilde{\nabla }_{%
\overline{j}}^{H}X_{\overline{i}}=0 \\ 
\end{array}%
\right.  \tag{13}
\end{equation}%
with respect to the adapted frame in $T\left( M_{n}\right) $.

Since we have

\begin{equation*}
\nabla _{i}\nabla _{l}X_{j}+\left( R_{hsij}+R_{hsji}\right) X^{h}=0
\end{equation*}%
as a consequence of $\nabla _{i}X_{j}+\nabla _{j}X_{i}=0$ $\left( see\cite{4}%
\right) $, we conclude by means of $\left( 16\right) $ that the complete
lift $^{C}X$ is a Killing vector field in $T\left( M_{n}\right) $ if only if 
$X$ is a Killing vector field in $M_{n}.$

We next have

\begin{equation*}
\left( R_{hsij}+R_{hsji}\right) X^{h}=0\ \ and\ \ \ \ \left(
R_{sij}^{h}+R_{sji}^{h}\right) X_{h}=0
\end{equation*}%
as a consequence of the vanishing of the second covariant derivative of $X$ .

Summing up these results, we have

\begin{theorem}
Necessary and sufficient condition in order that $\left( a\right) $
complete, $\left( b\right) $ horizontal lifts to $T\left( M_{n}\right) $
with the metric $II+III$ , of a vector field $X$ in $M_{n}$, be a Killing
vector field \ in $T\left( M_{n}\right) $ are respectively that, $\left(
a\right) $ $X$ is a Killing vector field with vanishing covaryant derivative
in $M_{n}$, $\left( b\right) $ $X$ is a Killing vector field with vanishing
second covaryant derivative in $M_{n}$.
\end{theorem}

\begin{center}
\bigskip

\bigskip

\bigskip
\end{center}

\bigskip

\bigskip

\bigskip

\bigskip

\bigskip

\end{document}